\documentclass[graybox]{svmult}
\usepackage{mathptmx}       
\usepackage{helvet}         
\usepackage{courier}        
\usepackage{type1cm}        
\usepackage{makeidx}         
\usepackage{graphicx}        
\usepackage{multicol}        
\usepackage[bottom]{footmisc}

\usepackage{latexsym}
\usepackage{amsfonts}
\usepackage{enumerate}
\usepackage{multicol}
\usepackage{graphicx}
\usepackage{amssymb}
\usepackage{amsmath}
\usepackage{epic}

\makeindex             
\newtheorem{conj}[theorem]{Conjecture}

\begin{document}

\title{Open Problems About Sumsets in Finite Abelian Groups: Minimum Sizes and Critical Numbers}
\titlerunning{Open Problems About Sumsets in Finite Abelian Groups} 
\author{B\'{e}la Bajnok}
\institute{B\'{e}la Bajnok \at Department of Mathematics, Gettysburg College, Gettysburg, PA 17325-1486 USA \\ \email{bbajnok@gettysburg.edu}}
%
%
\maketitle

\abstract*{For a positive integer $h$ and a subset $A$ of a given finite abelian group,
we let $hA$, $h \hat{\;} A$, and $h_{\pm}A$ denote the $h$-fold sumset, restricted sumset, and signed
sumset of $A$, respectively. Here we review some of what is known and not yet known
about the minimum sizes of these three types of sumsets, as well as their corresponding
critical numbers.  In particular, we discuss several new open direct and
inverse problems.}

\abstract{For a positive integer $h$ and a subset $A$ of a given finite abelian group,
we let $hA$, $h \hat{\;} A$, and $h_{\pm}A$ denote the $h$-fold sumset, restricted sumset, and signed
sumset of $A$, respectively. Here we review some of what is known and not yet known
about the minimum sizes of these three types of sumsets, as well as their corresponding
critical numbers.  In particular, we discuss several new open direct and
inverse problems.}

\section{Introduction and notations}

Throughout this paper, $G$ denotes a finite abelian group of order $n \geq 2$, written in additive notation.  If $G$ is cyclic, we identify it with $\mathbb{Z}_n=\mathbb{Z}/n \mathbb{Z}$. We say that $G$ has type $(n_1,\dots,n_r)$ if $$G \cong \mathbb{Z}_{n_1} \times \cdots \times \mathbb{Z}_{n_r}$$ for integers $2 \leq n_1 | n_2 | \cdots | n_r$; here $r$ is the rank and $n_r$ is the exponent of $G$.  For an $m$-subset $A=\{a_1, \dots, a_m\}$ of $G$ and for a nonnegative integer $h$ we consider three types of sumsets:
\begin{itemize}
  \item the {\em $h$-fold sumset}: $$hA=\left\{ \Sigma_{i=1}^m \lambda_i a_i \mid \lambda_1, \cdots, \lambda_m \in \mathbb{N}_0, \; \Sigma_{i=1}^m \lambda_i=h \right\},$$
\item the {\em $h$-fold restricted sumset}: $$h \hat{\;} A=\left\{ \Sigma_{i=1}^m \lambda_i a_i \mid \lambda_1, \cdots, \lambda_m \in \{0,1\}, \; \Sigma_{i=1}^m \lambda_i=h \right\},$$
\item the {\em $h$-fold signed sumset}: $$h_{\pm} A=\left\{ \Sigma_{i=1}^m \lambda_i a_i \mid \lambda_1, \cdots, \lambda_m \in \mathbb{Z}, \; \Sigma_{i=1}^m |\lambda_i|=h \right\}.$$
\end{itemize}
We denote the set formed by the inverses of the elements of $A$ by $-A$; we say that $A$ is {\em symmetric} if $A=-A$ and that $A$ is {\em asymmetric} if $A$ and $-A$ are disjoint.  For an element $b \in \mathbb{Z}$, we write $b \cdot A$ for the {\em $b$-fold dilation} $\{b \cdot a_1, \dots, b \cdot a_m\}$ of $A$.  The subgroup of $G$ generated by $A$ is denoted by  $\langle A \rangle$.

It is a central question in additive combinatorics to evaluate minimum sumset sizes, in particular, for given $G$, $h$, and $1 \leq m \leq n$ the quantities
\begin{eqnarray*}
\rho(G,m,h) & = & \min \{ |hA| \mid A \subseteq G, |A|=m \}, \\ \\
\rho \hat{\;} (G,m,h) & = & \min \{ |h \hat{\;} A| \mid A \subseteq G, |A|=m \}, \\ \\
\rho_{\pm} (G,m,h) & = & \min \{ |h_{\pm}A| \mid A \subseteq G, |A|=m \}.
\end{eqnarray*}
Trivially, each value is 1 whenever $h=0$, and each value equals $m$ whenever $h=1$.  (To see that $\rho_{\pm} (G,m,1)=m$ note that every group has a symmetric subset of any size $m \leq n$.)  Below we assume that $h \geq 2$; in the case of restricted sums, we may and will also assume that $h \leq m-2$.

The study of $\rho(G,m,h)$ goes back two hundred years to the work of Cauchy \cite{Cau:1813a} who determined it for groups of prime order and $h=2$, and is now known for all parameters---see Section \ref{min size sums} below.  However,  
only partial results have been found for  $\rho \hat{\;} (G,m,h)$ and $\rho_{\pm}  (G,m,h)$---we discuss these in Sections \ref{min size restricted} and \ref{min size signed}.

We also consider minimum sumset sizes without restrictions on the number of terms added: 
\begin{eqnarray*}
\rho(G,m,\mathbb{N}_0) & = & \min \{ | \cup_{h=0}^{\infty} hA| \mid A \subseteq G, |A|=m \}, \\ \\
\rho \hat{\;} (G,m,\mathbb{N}_0) & = & \min \{ |  \cup_{h=0}^{\infty} h \hat{\;} A| \mid A \subseteq G, |A|=m \}, \\ \\
\rho_{\pm} (G,m,\mathbb{N}_0) & = & \min \{ |  \cup_{h=0}^{\infty} h_{\pm}A| \mid A \subseteq G, |A|=m \}.
\end{eqnarray*}
Since $\cup_{h=0}^{\infty} hA$ and $\cup_{h=0}^{\infty} h_{\pm}A$ both equal $\langle A \rangle,$ 
we have $$\rho(G,m,\mathbb{N}_0)=\rho_{\pm} (G,m,\mathbb{N}_0)= \min \{ d \in D(n) \mid d \geq m\},$$ where $D(n)$ is the set of positive divisors of $n$.
The set $\cup_{h=0}^{\infty} h \hat{\;} A$, also denoted by $\Sigma A$, is less understood; we discuss $\rho \hat{\;} (G,m,\mathbb{N}_0)$ in Section \ref{min size any number of terms} below.

Related to each function, we study the corresponding {\em critical number}: the minimum value of $m$, if it exists, for which the corresponding sumset of any  $m$-subset of $G$ is $G$ itself.   The study of critical numbers originated with the 1964 paper \cite{ErdHei:1964a} of Erd\H{o}s and Heilbronn; in Section \ref{critical section} below we review what is known and not yet known about them.

Furthermore, we also examine the so-called {\em inverse problems} corresponding to some of these quantities, that is, we look for subsets of the group that achieve the extremal values of our functions.

The open problems mentioned here are just some of the many intriguing questions about sumsets.

\section{Minimum size of $h$-fold sumsets}  \label{min size sums}

Given $G$, $m$, and $h$, which $m$-subsets of $G$ have the smallest $h$-fold sumsets?  Two ideas come to mind: place the elements into a coset of some subgroup, or have the elements form an arithmetic progression.  We may also combine these two ideas; for example, in the cyclic group $\mathbb{Z}_n$,  we take an arithmetic progression of cosets, as follows.  

For any divisor $d$ of $n$, we take the subgroup $H=\cup_{j=0}^{d-1} \{j \cdot n/d\}$, then set
$$A_d(n,m)= \cup_{i=0}^{c-1} (i+H) \; \bigcup \; \cup_{j=0}^{k-1} \{c + j \cdot n/d\},$$ 
where $m=cd+k$ and $1 \leq k \leq d$.  An easy computation shows that
$$|hA_d(n,m)|=\min\{n, \; (hc+1)d, \; hm-h+1\};$$         
letting $$f_d(m,h)= (hc+1)d = (h \lceil m/d \rceil -h+1)d,$$ and noting that $f_n(m,h)=n$ and $f_1(m,h)=hm-h+1$, allows us to write
$$|hA_d(n,m)|=\min\{f_n(m,h), \; f_d(m,h), \; f_1 (m,h)\}.$$ Therefore, with
$$u(n,m,h)=\min\{ f_d(m,h) \mid d \in D(n)\}$$ we get $\rho(\mathbb{Z}_n,m,h) \leq u(n,m,h).$  It turns out that a similar construction works in any group and that we cannot do better:  
\begin{theorem} [Plagne \cite{Pla:2006a}] \label{thm plagne rho}
For every $G$, $m$, and $h$, we have
$\rho(G,m,h) = u(n,m,h).$
\end{theorem}
(Here $u(n,m,h)$ is a relative of the Hopf--Stiefel function used also in topology and bilinear algebra; see, for example, \cite{EliKer:2005a}, \cite{Kar:2006a}, \cite{Pla:2003a}, and \cite{Sha:1984a}.)

With $\rho(G,m,h)$ thus determined, let us turn to the inverse problem of classifying all $m$-subsets $A$ of $G$ for which $hA$ has minimum size $\rho(G,m,h)$.  The general question seems complicated. For example, while one can show that for a 6-subset $A$ of $\mathbb{Z}_{15}$ to have a 2-fold sumset of size $\rho(\mathbb{Z}_{15},6,2)=9$, $A$ must be the union of two cosets of the order 3 subgroup of $\mathbb{Z}_{15}$, there are three different possibilities for $\rho(\mathbb{Z}_{15},7,2)=13$: $A$ can be the union of two cosets of the order 3 subgroup plus one additional element, or a coset of the order 5 subgroup together with two more elements, or an arithmetic progression of length 7. 

We are able to say more for $m$ values that are not more than the smallest prime divisor $p$ of $n$.  Note that, as a special case of Theorem \ref{thm plagne rho}, when $m \leq p$, we get 
\begin{eqnarray} \label{p, hm-h+1 eq}
\rho(G,m,h) & = & \min\{p,hm-h+1\}.
\end{eqnarray} The case when $p$ is greater than $hm-h+1$ easily follows from \cite{Kem:1960a}:
\begin{theorem} [Kemperman \cite{Kem:1960a}] \label{Kemp thm inverse}
Let $p$ be the smallest prime divisor of $n$, and assume that $h \geq 2$ and $p > hm-h+1$.  Then for an $m$-subset $A$ of $G$ we have $|hA|=\rho(G,m,h) =hm-h+1$ if, and only if, $A$ is an arithmetic progression in $G$.
\end{theorem}
For the case when $p$ is less than $hm-h+1$, we propose:
\begin{conj} \label{conj p < hm-h+1 inverse}
Let $p$ be the smallest prime divisor of $n$, and assume that $m \leq p < hm-h+1$.  Then for an $m$-subset $A$ of $G$ we have $|hA|=\rho(G,m,h) =p$ if, and only if, $A$ is contained in a coset of some subgroup $H$ of $G$ with $|H|=p$.
\end{conj}
This leaves the case when $p=hm-h+1$, where arithmetic progressions of length $m$ and $m$-subsets in a coset of a subgroup of order $p$ are two of several possibilities.  It may be an interesting problem to classify all such subsets.

\section{Minimum size of $h$-fold restricted sumsets} \label{min size restricted}

While the value of $\rho \hat{\;} (G,m,h)$ is not even known for cyclic groups in general, as it turns out, we get an extremely close approximation for it by considering the sets $A_d(n,m) \subseteq \mathbb{Z}_n$ of Section \ref{min size sums} above.  A somewhat tedious computation shows that we get
$$|h \hat{\;} A_d(n,m)|= \left\{
\begin{array}{ll}
\min\{n, \; (hc+1)  d, \; hm-h^2+1\} & \mbox{if} \; h \leq \min\{k,d-1\}; \\ \\
\min\{n, \; hm-h^2+1 + \delta_d\} & \mbox{otherwise},
\end{array} \right.$$
where $\delta_d$ is an explicitly computed {\em correction term} (see \cite{Baj:2013a} for details).  Letting 
$$u \hat{\;} (n,m,h)= \min \{|h \hat{\;} A_d(n,m)| \mid d \in D(n) \},$$ we thus get $\rho \hat{\;} (\mathbb{Z}_n,m,h) \leq u \hat{\;} (n,m,h)$.  Since  $|h \hat{\;} A_d(n,m)|$ equals $\min\{n, hm-h^2+1\}$ for both $d=1$ and $d=n$, we always have
$$\rho \hat{\;} (\mathbb{Z}_n,m,h) \leq \min\{n, hm-h^2+1\}.$$  As is well known, equality holds for prime $n$:

\begin{theorem} [Dias Da Silva, Hamidoune \cite{DiaHam:1994a}; Alon, Nathanson, Ruzsa \cite{AloNatRuz:1995a, AloNatRuz:1996a}] \label{Dias Da Silva and Hamidoune 1}
For a prime $p$ we have $$\rho\hat{\;} (\mathbb{Z}_p,m,h) = \mathrm{min} \{p, hm-h^2+1\}.$$

\end{theorem}  

The lower bound $u \hat{\;} (n,m,h)$ is surprisingly accurate for cyclic groups of composite order as well: For all $(n,m,h)$ with $n \leq 40$, we find that equality holds in over 99.9\% of cases, and when it does not, then $\rho \hat{\;} (\mathbb{Z}_n,m,h)$ and $u \hat{\;} (n,m,h)$ differ only by 1.  All the exceptions that are known come from the construction that we explain next. 

Recall that the $m$ elements in $A_d(n,m)$ are within $c+1=\lceil m/d \rceil $ cosets of the order $d$ subgroup $H$ of $\mathbb{Z}_n$, and at most one of these cosets is not contained entirely in $A_d(n,m)$.  We now consider the variation when the $m$ elements are still within $c+1$ cosets of $H$, but exactly two of the cosets do not lie entirely in our set.  In order to do so, we write $$m=k_1+(c-1)d+k_2$$ with positive integers $k_1$ and $k_2$; we assume that $k_1 < d$, $k_2 < d$, but $k_1+k_2 >d$.  We then set
$$B_d(n,m)=\cup_{j=0}^{k_1-1} \{j \cdot n/d\} \; \bigcup\;  \cup_{i=1}^{c-1} (i \cdot g +H) \; \bigcup \; \cup_{j=0}^{k_2-1} \{c \cdot g + (j_0+j) \cdot n/d\},$$
where $0 \leq j_0 \leq d-1$ and $g \in \mathbb{Z}_n$.  As it turns out, $|h \hat{\;} B_d(n,m)|$ is less than $|h \hat{\;} A_d(n,m)|$ in just three specific cases: when $h=2$, $n$ is divisible by $2m-2$, and $m-1$ is not a power of 2; when $h=3$, $m=6$, and $n$ is divisible by 10; and when $h$ is odd, $n$ is divisible by $hm-h^2$, and $m+2$ is divisible by $h+2$ (see \cite{Baj:2013a}).  Moreover, every known instance when $\rho \hat{\;} (\mathbb{Z}_n,mh)$ is less than $u \hat{\;} (n,m,h)$ arises as one of these three cases.  Letting
$$w \hat{\;} (n,m,h)= \min \{|h \hat{\;} B_d(n,m)| \mid d \in D(n) \},$$ we 
see that $\rho \hat{\;} (\mathbb{Z}_n,m,h)$ is at most $\min\{u \hat{\;} (n,m,h), w \hat{\;} (n,m,h)\},$  but we also believe that equality holds: 
\begin{conj} \label{conj rho hat u w}
For all $n$, $m$, and $h$, we have $$\rho \hat{\;} (\mathbb{Z}_n,m,h)= \min\{u \hat{\;} (n,m,h), w \hat{\;} (n,m,h)\}.$$  
\end{conj}

Let us highlight the case $h=2$.  First, note that Conjecture \ref{conj rho hat u w} then becomes:
$$\rho \hat{\;} (\mathbb{Z}_n,m,2)=\left\{
\begin{array}{ll}
\min\{\rho (\mathbb{Z}_n,m,2), 2m-4\} & \mbox{if} \; 2|n \; \mbox{and} \; 2|m, \\
& \mbox{or} \; (2m-2)|n \; \mbox{and} \; \log_2 (m-1) \not \in \mathbb{N}; \\ \\
\min\{\rho (\mathbb{Z}_n,m,2), 2m-3\} & \mbox{otherwise.}
\end{array} \right.$$
Some general inequalities are known: Plagne \cite{Pla:2006b} proved that the upper bound 
$$\rho \hat{\;} (G,m,2) \leq \min\{\rho (G,m,2), 2m-2\}$$ holds for all groups, and Eliahou and Kervaire \cite{EliKer:1998a} proved that the lower bound 
$$\rho \hat{\;} (G,m,2) \geq \min\{\rho (G,m,2), 2m-3\}$$ holds for all elementary abelian $p$-groups for odd $p$.  Furthermore, Lev \cite{Lev:2000a} conjectured the lower bound 
$$\rho \hat{\;} (G,m,2) \geq \min\{\rho (G,m,2), 2m-3-|\mathrm{Ord}(G,2)|\},$$ where $\mathrm{Ord}(G,2)$ is the set of elements of $G$ that have order 2, and Plagne \cite{Pla:2006b} conjectured that $\rho \hat{\;} (G,m,2)$ and $\rho (G,m,2)$ can differ by at most 2.  (We should add that no such statement is possible for higher $h$ values: as was proven in \cite{Baj:2013a}, when $h \geq 3$, for any $C \in \mathbb{N}$, one can find a group $G$ and a positive integer $m$ so that
$\rho \hat{\;} (G,m,h)$ and $\rho (G,m,h)$ differ by  $C$ or more.)    
 
As in Section \ref{min size sums}, we are able to say more when $m \leq p$ with $p$ being the smallest prime divisor of $n$.  We believe that the following analogue of (\ref{p, hm-h+1 eq}) holds:
\begin{conj} \label{conj p hm-h^2+1 eq}
If $p$ is the smallest prime divisor of $n$ and $h< m \leq p$, then 
$$\rho \hat{\;} (G,m,h)=\min\{p, hm-h^2+1\}.$$
\end{conj}
Note that Conjecture \ref{conj p hm-h^2+1 eq} is a generalization of Theorem \ref{Dias Da Silva and Hamidoune 1}.  

Turning to inverse problems: our analogues for Theorem \ref{Kemp thm inverse} and Conjecture \ref{conj p < hm-h+1 inverse} are:
\begin{conj} \label{conj p > hm-h^2+1}
Let $p$ be the smallest prime divisor of $n$, and assume that $2 \leq h \leq m-2$ and $p > hm-h^2+1$.  Then for an $m$-subset $A$ of $G$ we have $|h \hat{\;} A|=hm-h^2+1$ if, and only if, $h=2$, $m=4$, and $A=\{a,a+g_1,a+g_2,a+g_1+g_2\}$ for some $a, g_1, g_2 \in G$, or $A$ is an arithmetic progression in $G$.
\end{conj}

\begin{conj}
Let $p$ be the smallest prime divisor of $n$, and assume that $m \leq p < hm-h^2+1$.  Then for an $m$-subset $A$ of $G$ we have $|h \hat{\;} A| =p$ if, and only if, $A$ is contained in a coset of some subgroup $H$ of $G$ with $|H|=p$.
\end{conj}

K\'arolyi proved Conjectures \ref{conj p hm-h^2+1 eq} and \ref{conj p > hm-h^2+1} for $h=2$ \cite{Kar:2003a,Kar:2004a}.

\section{Minimum size of $h$-fold signed sumsets} \label{min size signed}

Studying the function $\rho_{\pm}(G,m,h)$ provides us with several surprises.  First, we realize that, unlike it is the case for $\rho(G,m,h)$, the value of $\rho_{\pm}(G,m,h)$ depends on the structure of $G$ and not just on the order $n$ of $G$.   Second, while the size of the signed sumset of a subset is usually much greater than the size of its sumset, the value of $\rho_{\pm}(G,m,h)$ equals $\rho(G,m,h)$ surprisingly often; in fact, there is only one case with $n \leq 24$ where the two are not equal: $\rho_{\pm}(\mathbb{Z}_3^2,4,2)=8$ while $\rho(\mathbb{Z}_3^2,4,2)=7$.  Furthermore, one might think that symmetric sets provide the smallest minimum size, but sometimes asymmetric sets or even {\em near-symmetric sets}---sets that become symmetric by the removal of one element---are better; we are able to prove, though, that one of these three types always provides the minimum size.    

For our treatment below, we use the functions $f_d(m,h)$ and $u(n,m,h)$ defined in Section \ref{min size sums}.  For cyclic groups, we have the following result:
\begin{theorem} [Bajnok and Matzke \cite{BajMat:2014a}] \label{pm for cyclic} 
For cyclic groups $G$, $m$, and $h$, we have $\rho_{\pm}(G,m,h)=\rho(G,m,h)$.

\end{theorem}
The proof of Theorem \ref{pm for cyclic} follows from the fact that for each $d \in D(n)$ one can find a symmetric subset $R$ of $G$ of size at least (but not necessarily equal to) $m$ for which $|hR| \leq f_d(n,m)$.  

More generally, for a group of type $(n_1,\dots,n_r)$ one can prove that 
$$\rho_{\pm}(G,m,h) \leq \min \{ \Pi_{i=1}^r \rho_{\pm}(\mathbb{Z}_{n_i},m_i,h) \mid m_i \leq n_i, \Pi_{i=1}^r m_i \geq m \},$$ so by Theorems \ref{thm plagne rho} and 
 \ref{pm for cyclic}, 
$$\rho_{\pm}(G,m,h) \leq \min \{ \Pi_{i=1}^r u(n_i, m_i, h) \mid m_i \leq n_i, \Pi_{i=1}^r m_i \geq m \}.$$
Furthermore, in \cite{BajMat:2014a} we proved that
$$\min \{ \Pi_{i=1}^r u(n_i, m_i, h) \mid m_i \leq n_i, \Pi_{i=1}^r m_i \geq m \}=\min\{f_d(m,h) \mid d \in D(G,m)\},$$
 where $D(G,m)$ consists of all $d \in D(n)$ that can be written as $d=\Pi_{i=1}^r d_i$ with $d_i \in D(n_i)$ and $dn_r \geq d_rm$.  (We may observe that for cyclic groups $D(G,m)=D(n)$.)  Letting
$$u_{\pm} (G,m,h)=\min\{f_d(m,h) \mid d \in D(G,m)\}$$ thus results in the upper bound
$\rho_{\pm}(G,m,h) \leq u_{\pm} (G,m,h)$.  Of course, we also have $\rho_{\pm}(G,m,h) \geq u  (n,m,h)$, so to get lower and upper bounds for $\rho_{\pm}(G,m,h)$, one can minimize the values of $f_d(m,h)$ for all $d \in D(n)$ and for all $d \in D(G,m)$, respectively.  In fact, with one specific exception, that we are about to explain, we believe that $\rho_{\pm}(G,m,h) = u_{\pm}  (G,m,h)$  holds for all $G$, $m$, and $h$.

We can observe that if $A$ is asymmetric, then $0 \not \in 2_{\pm}A$.  Consequently, if $d$ is an odd divisor of $n$ and $d \geq 2m+1$, then we can choose an $m$-subset of $G$ whose 2-fold signed sumset has size less than $d$, and thus $\rho_{\pm}(G,m,2) \leq d-1$. We believe that this is the only possibility for $\rho_{\pm}(G,m,h)$ to be less than $u_{\pm}  (G,m,h)$:
\begin{conj} [Bajnok and Matzke \cite{BajMat:2014a}] \label{pm conj value}
For all $G$, $m$, and $h \geq 3$, we have $\rho_{\pm}(G,m,h) = u_{\pm}  (G,m,h)$.  

Furthermore, with $D_o(n)$ denoting the set of odd divisors of $n$ that are greater than $2m$, we have
$$\rho_{\pm}(G,m,2) = \left\{
\begin{array}{ll}
u_{\pm}  (G,m,2) & \mbox{if} \; D_o(n) = \emptyset, \\ \\
\min\{u_{\pm}  (G,m,2), d_m-1\} & \mbox{if} \; d_m=\min D_o(n).
\end{array}
\right.$$

  \end{conj}

We can say more about elementary abelian groups.  Clearly, $\rho_{\pm}(\mathbb{Z}_2^r,m,h)=\rho(\mathbb{Z}_2^r,m,h)$, so consider $\mathbb{Z}_p^r$ where  $p$ is an odd prime.  When $p \leq h$, one can prove that $\rho_{\pm}(\mathbb{Z}_p^r,m,h)=\rho(\mathbb{Z}_p^r,m,h)$ \cite{BajMat:2014b}.  The case when $h$ is less than $p$ is more delicate; we need the following notations.  First, set $k$ equal to the largest integer for which $p^k+\delta \leq hm-h+1$, where $\delta=0$ if $p-1$ is divisible by $h$ and $\delta=1$ otherwise.  Second, set $q$ equal to the largest integer  for which $(hq+1)p^k+\delta \leq hm-h+1$.  With these notations, we have the following result:
\begin{theorem} [Bajnok and Matzke \cite{BajMat:2014b}]  Suppose that either $p \leq h$, or that $h<p$ and $m \leq (q+1)p^k$
with $k$ and $q$ defined as above.  Then $\rho_{\pm}(\mathbb{Z}_p^r,m,h)=\rho(\mathbb{Z}_p^r,m,h)$. 
\end{theorem} 
We believe that $\rho_{\pm}(\mathbb{Z}_p^r,m,h)$ is greater than $\rho(\mathbb{Z}_p^r,m,h)$ in the remaining case:
\begin{conj} [Bajnok and Matzke \cite{BajMat:2014b}]  \label{BajMat conj for p groups}
If $h<p$ and $m > (q+1)p^k$
with $k$ and $q$ defined as above, then $\rho_{\pm}(\mathbb{Z}_p^r,m,h) > \rho(\mathbb{Z}_p^r,m,h)$. 
\end{conj} 

Using Vosper's Theorem \cite{Vos:1956a} and (Lev's improvement \cite{Lev:2006a} of) Kemperman's  results on so-called {\em critical pairs} \cite{Kem:1960a}, in \cite{BajMat:2014b} we were able to prove Conjecture \ref{BajMat conj for p groups} for the case when $r=2$ and $h=2$; therefore we have a complete account for all $m$ for which $\rho_{\pm}(\mathbb{Z}_p^2,m,2) = \rho(\mathbb{Z}_p^2,m,2)$.  In particular, we found that there are exactly $(p-1)^2/4$ values of $m$ where equality does not hold.  We have not been able to find any groups where this proportion is higher than $1/4$, and believe that there are none:

\begin{conj}
For any abelian group of order $n$, $\rho_{\pm} (G, m, 2)$ and $\rho(G, m, 2)$ disagree for fewer than $n/4$ values of $m$.

\end{conj} 

Let us turn now to the inverse problem of classifying all $m$-subsets $A$ of $G$ for which $|h_{\pm}A|=\rho_{\pm}(G,m,h)$.  Letting  $\mathrm{Sym}(G,m)$, $\mathrm{Nsym}(G,m)$, and $\mathrm{Asym}(G,m)$ denote the collection of $m$-subsets of $G$ that are, respectively, symmetric, near-symmetric (that is, become symmetric after removing one element), and asymmetric, in \cite{BajMat:2014a} we proved that 
$$\rho_{\pm}(G,m,h)= \min \{|h_{\pm}A| \mid A \in \mathrm{Sym}(G,m) \cup \mathrm{Nsym}(G,m) \cup \mathrm{Asym}(G,m) \}.$$  (None of the three types are superfluous.)  This does not completely solve the inverse problem: we may have other subsets with $|h_{\pm}A|=\rho_{\pm}(G,m,h)$.  Furthermore, it would be interesting to know exactly when each of the three types of sets just described yields a signed sumset of minimum size.

\section{Minimum size of restricted sumsets with an arbitrary number of terms} \label{min size any number of terms}

In this section we attempt to find the minimum size $\rho \hat{\;} (G,m,\mathbb{N}_0)$ of $\Sigma A=\cup_{h=0}^{\infty} h \hat{\;} A$ among all $m$-subsets of $G$.  We restrict our attention to cyclic groups.

As before, we choose a divisor $d$ of $n$, and consider an arithmetic progression of cosets of the subgroup $H$ of order $d$ in $\mathbb{Z}_n$.  We again write $m=cd+k$ with $1 \leq k \leq d$, and construct a set $C_d(n,m)$ that lies in exactly $c+1 = \lceil m/d \rceil $ cosets of $H$, as follows.

Assume first that $c$ is even.  In this case, we let $C_d(n,m)$ consist of the collection of $c$ cosets 
$$\{i+H \mid -c/2 \leq i \leq c/2-1\},$$ together with $k$ elements of the coset $c/2+H$.  (It makes no difference which $k$ elements we choose.) 
It is easy to see then that 
$$\Sigma C_d(n,m) = \{i+H \mid -(c^2+2c)/8 \cdot d \leq i \leq (c^2-2c)/8 \cdot d  +c/2 \cdot k\},$$
and thus 
\begin{eqnarray*}
|\Sigma C_d(n,m)| & = & \min \left\{n, \; (c^2/4 \cdot d + c/2 \cdot k+1 ) \cdot d \right\}\\
& = & \min \left\{n, \; (c/2 \cdot m - c^2/4 \cdot d + 1 ) \cdot d \right\}.
\end{eqnarray*}

Similarly, when $c$ is odd, we set $C_d(n,m)$ equal to the collection
$$\{i+H \mid -(c-1)/2 \leq i \leq (c-1)/2\},$$ together with $k$ elements of the coset $(c+1)/2+H$; this time we find that
$$|\Sigma C_d(n,m)|=\min \left\{n, \; \left((c+1)/2 \cdot m - (c+1)^2/4 \cdot d + 1 \right) \cdot d \right\}.$$
Therefore, letting
\begin{eqnarray*}
F_d(m) & = & \left( \lceil c/2 \rceil \cdot m - \lceil c/2 \rceil^2 \cdot d +1 \right) \cdot d \\
& = & \left( \left \lceil (m/d-1)/2 \right \rceil \cdot m - 
\left \lceil (m/d-1)/2  \right \rceil ^2 \cdot d + 1 \right) \cdot d,
\end{eqnarray*}
and noting that $F_n(m)=n$, we may write $|\Sigma C_d(n,m)|=\min \{F_n(m), F_d(m)\}$.  Setting $$u(n,m,\mathbb{N}_0)=\min \{F_d(m) \mid d \in D(n)\},$$ we get:
\begin{theorem} \label{sigma lower cyclic 1}
For all positive integers $n$ and $m \leq n$, we have $\rho \hat{\;} (\mathbb{Z}_n,m,\mathbb{N}_0) \leq u(n,m,\mathbb{N}_0).$
\end{theorem}
After some numerical experimentation, we believe that equality holds:
\begin{conj} \label{conj sigma cyclic 1}
For all positive integers $n$ and $m \leq n$, we have $\rho \hat{\;} (\mathbb{Z}_n,m,\mathbb{N}_0) = u(n,m,\mathbb{N}_0).$
\end{conj}
Note that, since $F_1(m)=\lfloor m^2/4 \rfloor ^2+1,$ Theorem \ref{sigma lower cyclic 1} implies that 
\begin{eqnarray} \label{F_1, F_n lower}
\rho \hat{\;} (\mathbb{Z}_n,m,\mathbb{N}_0) & \leq & \min \left\{ n, \; \lfloor m^2/4 \rfloor +1 \right\}
\end{eqnarray} holds for all $n$ and $m \leq n$.

Next, we examine groups of prime order.  Let $p$ be a positive prime.  Trivially, for any subset $A$ and any positive integer $h$, the $h$-fold restricted sumset of $A$ is contained in $\Sigma A$ and, therefore,   
$\rho \hat{\;} (\mathbb{Z}_p,m,\mathbb{N}_0)$ cannot be less than $\rho \hat{\;} (\mathbb{Z}_p,m, \lfloor m/2 \rfloor)$.  By Theorem \ref{Dias Da Silva and Hamidoune 1}, 
$$\rho \hat{\;} (\mathbb{Z}_p,m, \lfloor m/2 \rfloor) = \min \left \{ p, \; \lfloor m/2 \rfloor \cdot m - \lfloor m/2 \rfloor^2 +1 \right\}= \min \left \{ p, \; \lfloor m^2/4 \rfloor  +1 \right\};$$
together with our upper bound (\ref{F_1, F_n lower}), we arrive at:
\begin{theorem}
Conjecture \ref{conj sigma cyclic 1} holds for groups of prime order; in particular, $$\rho \hat{\;} (\mathbb{Z}_p,m, \mathbb{N}_0)= \min \left \{ p, \; \lfloor m^2/4 \rfloor  +1 \right\}$$ for all primes $p$ and $m \leq p$.
\end{theorem}

There have been some studies of several variations of $\rho \hat{\;} (G,m,\mathbb{N}_0)$ provided by various restrictions on the subsets $A$ of $G$.  We mention only one pair of such results. Recall that $\mathrm{Asym}(G,m)$ denotes the collection of asymmetrical $m$-subsets of $G$; also set $\Sigma ^*A=\cup_{h=1}^{\infty} h \hat{\;} A$.  Furthermore, let
$$\rho_{\mathrm{A}} \hat{\;} (G,m,\mathbb{N}_0)=\min \{ | \Sigma A | \mid A \in \mathrm{Asym}(G,m)\},$$
$$\rho_{\mathrm{A}} \hat{\;} (G,m,\mathbb{N})=\min \{ | \Sigma^* A | \mid A \in \mathrm{Asym}(G,m)\}.$$
With these notations:
\begin{theorem} [Balandraud \cite{Bal:2012a}, \cite{Bal:2012b}, \cite{Bal:2017a}]
For every odd prime $p$ and every $m \leq (p-1)/2$ we have 
$$\rho_{\mathrm{A}} \hat{\;} (\mathbb{Z}_p,m,\mathbb{N}_0)=\min \{ p, (m^2+m)/2+1\},$$
$$\rho_{\mathrm{A}} \hat{\;} (\mathbb{Z}_p,m,\mathbb{N})=\min \{ p, (m^2+m)/2 \}.$$
\end{theorem}  
The fact that the  values are upper bounds is provided by the set $\{1,2,\dots,m\}$.

\section{Critical numbers}  \label{critical section}

Given $G$ and $h$, we define the {\em $h$-critical number} $\chi (G,h)$ as the least integer $m$ for which $hA=G$ holds for all $m$-subsets $A$ of $G$; we define $\chi \hat{\;}(G,h)$ and $\chi_{\pm} (G,h)$ analogously.  We also define the {\em critical number} $\chi (G, \mathbb{N}_0)$ as the smallest value of $m$ for which $\cup_{h=0}^{\infty} hA=G$ holds for all $m$-subsets $A$ of $G$; we define $\chi \hat{\;}(G, \mathbb{N}_0)$ and $\chi_{\pm} (G, \mathbb{N}_0)$ analogously.  

The study of critical numbers originated with the 1964 paper \cite{ErdHei:1964a} of Erd\H{o}s and Heilbronn: they studied the variation (in groups of prime order) where only $m$-subsets of $G \setminus \{0\}$ were considered.  (As we now know, the restriction to subsets that do not contain 0 does not change the critical numbers when the number of terms is a fixed value of $h$, but reduces them by 1 when the number of terms is arbitrary; for example, 
the least integer $m$ for which $hA=G$ holds for all $m$-subsets $A$ of $G \setminus \{0\}$ equals $\chi (G,h),$ but the least integer $m$ for which $ \Sigma A=G$ holds for all $m$-subsets $A$ of $G \setminus \{0\}$ equals $\chi \hat{\;}(G,\mathbb{N}_0)-1$; see \cite{Baj:2015a}.)

Two of these six quantities are obvious: Since $\cup_{h=0}^{\infty} hA$ and $\cup_{h=0}^{\infty} h_{\pm}A$ are both equal to $\langle A \rangle$, we have
$$\chi (G, \mathbb{N}_0)=\chi_{\pm} (G, \mathbb{N}_0)=n/p+1,$$ where $p$ is the smallest prime divisor of $n$.  Furthermore, $\chi (G,h)$ and $\chi \hat{\;} (G,\mathbb{N}_0)$ have now been determined, but the remaining two quantities are not known in general.  Let us review what we know. 

To state the result for $\chi (G,h)$, we need to introduce the---perhaps already familiar---function  \label{vdef}
$$v_g(n,h)= \max \left\{ \left( \left \lfloor \left(d-1-\mathrm{gcd} (d, g) \right)/h \right \rfloor +1  \right) \cdot n/d  \; : \; d \in D(n) \right\}$$ ($n,g,h \in \mathbb{N}$).    
We should note that this function has appeared elsewhere in additive combinatorics already.  For example, according to the classical result of 
Diananda and Yap \cite{DiaYap:1969a}, the maximum size of a sum-free set (that is, a set $A$ that is disjoint from $2A$) in the cyclic group $\mathbb{Z}_n$ is given by
$$v_1(n,3) =\left\{
\begin{array}{ll}
\left(1+1/p \right) \cdot n/3 & \mbox{if $n$ has prime divisors $p \equiv 2$ mod 3} \\ & \mbox{and $p$ is the smallest,}\\ \\
\left\lfloor n/3 \right\rfloor & \mbox{otherwise.}\\
\end{array}\right.$$
Similarly, we proved in \cite{Baj:2009a} that the maximum size of a $(3,1)$-sum-free set in $\mathbb{Z}_n$ (where $A$ is disjoint from $3A$) 
equals $v_2(n,4)$.  More generally, $v_{k-l}(n,k+l)$ provides a lower bound for the maximum size of $(k,l)$-sum-free sets in $\mathbb{Z}_n$ (where $kA \cap lA=\emptyset$ for positive integers $k>l$) (see \cite{Baj:2009a}); equality holds in the case when $k-l$ and $n$ are relatively prime (see the paper \cite{HamPla:2003a} of Hamidoune and Plagne).  We can now state our result for $\chi (G,h)$:
\begin{theorem} [Bajnok \cite{Baj:2015a}]
For all $G$ and $h$, we have $\chi (G,h)=v_1(n,h)+1.$

\end{theorem}

Let us now see what we can say about $\chi \hat{\;} (G,h)$.  First, we can prove that ${\chi} \hat{\;} (G, h)$ is  well defined, except when $h \in \{2,n-2\}$ and $G$ is isomorphic to an elementary abelian 2-group.  Furthermore, for all $G$ with exponent at least 3, we have
$$\chi \hat{\;} (G,2)= (n+|\mathrm{Ord}(G,2)|+1)/2+1$$ and, as a consequence, when $h \geq (n+|\mathrm{Ord}(G,2)|-1)/2$, we have $\chi \hat{\;} (G,h) = h+2$ \cite{Baj:2015a}.  Regarding other values of $h$, few exact results are known; in particular, for $3 \leq h \leq \lfloor n/2 \rfloor-1$, we only know the value of $\chi \hat{\;} (\mathbb{Z}_n,h)$  when $n$ is prime or even.  

Indeed, for prime values of $p$, Theorem \ref{Dias Da Silva and Hamidoune 1} allows us to derive that 
$$\chi \hat{\;} (\mathbb{Z}_p,h)=\left \lfloor (p-2)/h \right \rfloor +h+1.$$
The case of even $n$ and $h=3$ was established by Gallardo, Grekos, et al.~in ~\cite{GalGre:2002a}; we generalized this in \cite{Baj:2015a}  (see also \cite{Baj:2017a}) to prove that
that for any $h$ and even $n \geq 12$, we have
$$\chi \hat{\;} (\mathbb{Z}_{n},h) = \left \{
\begin{array}{cl}
n/2+1 & \mbox{if} \; 3 \leq h \leq  n/2-2; \\ \\
n/2+2 & \mbox{if} \; h=n/2-1.
\end{array}
\right.$$

Let us take a closer look at the case of $h=3$.  In \cite{Baj:2015a} we proved that if $n \geq 16$ and has prime divisors congruent to $2$ mod $3$ and $p$ is the smallest such divisor, then 
  $$\chi \hat{\;} (\mathbb{Z}_n,3) \geq 
\left\{
\begin{array}{ll}
\left(1+1/p \right) \cdot n/3 +3 & \mbox{if} \; n=p,  \\ \\
\left(1+1/p \right) \cdot n/3 +2 & \mbox{if} \; n=3p,  \\ \\
\left(1+1/p \right) \cdot n/3 +1 & \mbox{otherwise}; 
\end{array} \right.$$
and if $n$ has no prime divisors congruent to $2$ mod $3$, then
  $$\chi \hat{\;} (\mathbb{Z}_n,3) \geq 
\left\{
\begin{array}{ll}
\left \lfloor n/3 \right \rfloor +4 & \mbox{if $n$ is divisible by $9$},  \\ \\
\left \lfloor n/3 \right \rfloor +3 & \mbox{otherwise}. 
\end{array} \right.$$
We also believe that, actually, equality holds above for all $n$---this is certainly the case if $n$ is even or prime; we have verified this (by computer) for all $n \leq 50$.  Our conjecture is a generalization of the one made by Gallardo, Grekos, et al.~in \cite{GalGre:2002a} that was proved (for large $n$) by Lev in \cite{Lev:2002a}.

The study of ${\chi} \hat{\;} (G, \mathbb{N}_0)$ posed considerable amount of challenges, but after several decades of attempts, due to the combined results of Diderrich and Mann \cite{DidMan:1973a}, Diderrich \cite{Did:1975a},  Mann and Wou \cite{ManWou:1986a}, Dias Da Silva and Hamidoune \cite{DiaHam:1994a}, Gao and Hamidoune \cite{GaoHam:1999a}, Griggs \cite{Gri:2001a}, and  Freeze, Gao, and Geroldinger \cite{FreGaoGer:2009a, FreGaoGer:2015a}, we have the value for every group:

\begin{theorem} [The combined results of authors above] \label{combined critical}

Suppose that $n \geq 10$, and let $p$ be the smallest prime divisor of $n$.  Then
$${\chi} \hat{\;} (G, \mathbb{N}_0)=\left\{
\begin{array}{ll}
\lfloor 2 \sqrt{n-2} \rfloor+1  & \mbox{if $G$ is cyclic of order $n=p$ or $n=pq$ where} \\
& \mbox{$q$ is prime and $3 \leq p \leq q \leq p+\lfloor 2 \sqrt{p-2} \rfloor+1$}\footnotemark, \\  \\
n/p+p-1 & \mbox{otherwise}.
\end{array}
\right.
$$

\end{theorem}
\footnotetext{Note that $\lfloor 2 \sqrt{n-2} \rfloor+1=n/p+p$ in this case.}

In closing, we state an intriguing question for the inverse problem regarding ${\chi} \hat{\;} (\mathbb{Z}_p, \mathbb{N}_0)$, that is, the attempt to classify all subsets $A$ of size 
${\chi} \hat{\;} (\mathbb{Z}_p, \mathbb{N}_0)-1=\lfloor 2 \sqrt{p-2} \rfloor$ in the cyclic group $\mathbb{Z}_p$ of odd prime order $p$ for which $\Sigma A \neq \mathbb{Z}_p$.  First, some notations and an observation.  Following standard techniques, we identify the elements of $\mathbb{Z}_p$ with integers between $-(p-1)/2$ and $(p-1)/2$ (inclusive), therefore we can write $A=A_1 \cup A_2$ where $A_1$ consists of the nonnegative elements of $A$, and $A_2$ consists of its negative elements. We define the {\em norm} of $A \subseteq \mathbb{Z}_p$, denoted by $||A||$, as the sum of the absolute values of its elements, thus $$||A||=\Sigma_{a \in A_1} a-\Sigma_{a \in A_2}a.$$  We note that if $||A|| \leq p-2$, then $$1+\Sigma_{a \in A_1} a \not \in \Sigma A;$$ in particular, $\Sigma A \neq \mathbb{Z}_p$.  Consequently, if  $||A|| \leq p-2$, then $\Sigma (b \cdot A) \neq \mathbb{Z}_p$ for any $b \in \mathbb{Z}_p$.  We believe that this simple condition answers our inverse problem for all large enough primes; namely: there is a positive integer $p_0$ so that if $p > p_0$ is prime and $A \subseteq \mathbb{Z}_p$ has size ${\chi} \hat{\;} (\mathbb{Z}_p, \mathbb{N}_0)-1=\lfloor 2 \sqrt{p-2} \rfloor$, then 
$\Sigma A \neq \mathbb{Z}_p$ if, and only if, there is a nonzero element $b \in \mathbb{Z}_p$ for which $|| b \cdot A|| \leq p-2$. 
(We verified that all primes under 40, with the exception of $p=17$, satisfy this condition.)    
We mention the following related result:

\begin{theorem} [Nguyen, Szemer\'edi, Vu \cite{NguSzeVu:2008a}]
Let $p$ be an odd prime, and let $A \subseteq \mathbb{Z}_p$ have size $|A| \geq 1.99 \sqrt{p}$.  If $\Sigma A \neq \mathbb{Z}_p$, then 
there is a nonzero element $b \in \mathbb{Z}_p$ for which $|| b \cdot A|| \leq p+ O(\sqrt{p})$.
\end{theorem}

\end{document}